\documentclass[11pt]{amsart}

\usepackage{amsmath, amssymb, amsthm}

\date{\today}

\newcommand{\la}{\langle}
\newcommand{\ra}{\rangle}

\theoremstyle{plain}
\newtheorem{theorem}{Theorem}[section]
\newtheorem{cor}[theorem]{Corollary}

\newtheorem{lemma}[theorem]{Lemma}
\newtheorem{remark}[theorem]{Remark}
\newtheorem{ex}[theorem]{Example}

\date{\today}
\title{On excesses of frames}
\author[D. Baki\' c and T. Beri\' c]
{Damir Baki\' c and Tomislav Beri\' c}

\address{Department of Mathematics, University of Zagreb,
Bijeni\v cka cesta 30, 10000 Zagreb, Croatia.}
\email{bakic@math.hr}
\email{tberic@math.hr}

\begin{document}

\begin{abstract}
We show that any two frames in a separable Hilbert space that are dual to each other have the same excess. Some new relations for the analysis resp.~synthesis operators of dual frames are also derived. We then prove that pseudo-dual frames and, in particular, approximately dual frames have the same excess. We also discuss various results on frames in which excesses of frames play an important role.

\vspace{.1in}

\noindent
{\it AMS Mathematics Subject Classification:} 42C15.

\noindent
{\it Key words and phrases:}  frame, Parseval frame, excess.
\end{abstract}

\maketitle

\section{Introduction and preliminaries}

Let $H$ be a separable Hilbert space with the inner product $\langle \cdot, \cdot \rangle$.
A sequence $(f_n)_{n=1}^{\infty}$ in $H$ is a frame if there exist positive constants $A$ and $B$, that are called frame bounds, such that
\begin{equation}\label{frame def}
A\|x\|^2\leq \sum_{n=1}^{\infty}|\la x,f_n\ra |^2\leq B\|x\|^2,\,\forall x \in H.
\end{equation}
Frame bounds are not unique. The optimal upper frame bound is the infimum over all upper frame bounds, and the optimal lower frame bound is the supremum over all lower frame bounds.  If $A=B$ we say that the frame is tight and, in particular, if $A=B=1$ so that
\begin{equation}
\sum_{n=1}^{\infty}|\la x,f_n\ra |^2=\|x\|^2,\,\forall x \in H,
\end{equation}
we say that $(f_n)_{n=1}^{\infty}$ is a Parseval frame.

Frames were first introduced by Duffin and Schaeffer (\cite{DS}). Today frames play important roles in many applications in mathematics, science and engineering. For general information and basic facts about frames we refer the reader to \cite{Cas}, \cite{Chr}, \cite{Dau} and \cite{HW}.

If only the second inequality in (\ref{frame def}) is satisfied, we say that $(f_n)_{n=1}^{\infty}$ is a Bessel sequence. For each Bessel sequence $(f_n)_{n=1}^{\infty}$ in $H$ one defines the analysis operator $U: H \rightarrow l^2$ by $Ux=(\langle x,f_n\rangle)_{n=1}^{\infty},\,x\in H$. It is evident that $U$ is bounded; moreover, if $(f_n)_{n=1}^{\infty}$ is a frame, then $U$ is also bounded from below. Its adjoint operator $U^*$, which is called the synthesis operator, is given by $U^*((c_n)_{n=1}^{\infty})=\sum_{n=1}^{\infty}c_nf_n,\,(c_n)_{n=1}^{\infty}\in l^2$. The synthesis operator $U^*$ is a surjection and the composition $U^*U$ (sometimes called the frame operator) is an invertible operator on $H$. It turns out that the sequence $((U^*U)^{-1}f_n)_{n=1}^{\infty}$ is also a frame for $H$ that satisfies
\begin{equation}\label{defdual}
x=\sum_{n=1}^{\infty}\la x,f_n\ra (U^*U)^{-1}f_n=\sum_{n=1}^{\infty}\la x,(U^*U)^{-1}f_n\ra f_n,\,\,\forall x \in H.
\end{equation}
In the light of this reconstruction formula we say that $((U^*U)^{-1}f_n)_{n=1}^{\infty}$ is the canonical dual frame of $(f_n)_{n=1}^{\infty}$.
The above equalities show that each frame is complete in $H$. However, a frame need not be a basis and the representations in (\ref{defdual}) need not be unique.
In general, a frame $(g_n)_{n=1}^{\infty}$ for $H$ is called a dual frame (or an alternate dual) for $(f_n)_{n=1}^{\infty}$ if we have
\begin{equation}\label{defdualgen}
x=\sum_{n=1}^{\infty}\la x,f_n\ra g_n=\sum_{n=1}^{\infty}\la x,g_n\ra f_n,\,\,\forall x \in H.
\end{equation}
If we denote by $U$ and $V$ the analysis operators of $(f_n)_{n=1}^{\infty}$ and $(g_n)_{n=1}^{\infty}$, respectively, then the above duality relations are simply described by the equality $V^*U=I$ (or, equivalently, $U^*V=I$).
Each frame $(f_n)_{n=1}^{\infty}$ for $H$ that is not a basis for $H$ has infinitely many dual frames.

 Frames that are not bases are overcomplete, {\it i.e.}, there exist their proper subsets which are complete. The excess $e((f_n)_{n=1}^{\infty})$ of the frame $(f_n)_{n=1}^{\infty}$ is defined in \cite{CHet_al} as the greatest integer $k$ such that $k$ elements can be deleted from the frame and still leave a complete set, or $+\infty$ if there is no upper bound to the number of elements that can be removed. By Lemma 4.1 from \cite{CHet_al} we have $e((f_n)_{n=1}^{\infty})=\mbox{dim}(\mbox{Ker}\,U^*)$, where $U$ is the analysis operator of $(f_n)_{n=1}^{\infty}$.

If $e((f_n)_{n=1}^{\infty})<\infty$ it can be shown (see \cite{Holub}) that the frame is simply a Riesz basis to which finitely many elements (in fact, precisely
$e((f_n)_{n=1}^{\infty})$ elements) have been adjoined. Such frames are called "near-Riesz bases".

\vspace{.1in}

The present paper is concerned with excesses of frames. In Theorem \ref{isti visak} we prove that dual frames have the same excess. Surprisingly, it seems that this fact has been overlooked by now, although its proof uses only basic tools. This also leads to some new relations for the analysis, resp.~synthesis operators of dual frames; see Corollaries \ref{dualsynthesis} and \ref{dualform}. We then prove that also pseudo-dual and approximately dual frames have the same excess. (Notions of pseudo-duals and approximate duals are explained in Section 2.) In addition, we include some results on frame perturbations, Parseval duals and the fundamental identity for Parseval frames (those that are known from the literature are presented with new proofs) in which excesses of frames play an important role.

\vspace{.1in}

It is worth noting here that the excess is a rather crude way of measuring the redundancy of frames. One can use a more refined measure as in~\cite{BL}, or restrict attention to a broad class of frames which are well behaved with regards to redundancy, e.g. localized frames (see~\cite{BCL}). The present paper, however, is concerned exclusively with excesses of frames.

\vspace{.1in}

We end this introductory section by fixing our notation. Throughout the paper $H$ denotes an infinite-dimensional Hilbert space. We shall also  tacitly assume that sequences in $H$ are indexed by the set $\Bbb N$ of natural numbers; so in the sequel we shall write $(f_n), (g_n)$ ... instead of $(f_n)_{n=1}^{\infty}$, $(g_n)_{n=1}^{\infty}$ ...

If $X$ and $Y$ are subspaces of $H$ with trivial intersection we denote by $X\stackrel{.}{+}Y$ their direct sum.

The identity operator is denoted by $I$. We denote by $\Bbb B(H,K)$ the space of all bounded operators of Hilbert spaces $H$ and $K$. The range and the kernel of an operator $T\in \Bbb B(H,K)$ are denoted by $\mbox{Im}\,T$ and $\mbox{Ker}\,T$, respectively.   If $T\in \Bbb B(H,K)$ has closed range, its pseudo-inverse is denoted by $T^{\dag}$. For basic facts concerning pseudo-inverses we refer to appendix A7 in \cite{Chr}. Finally, we denote by $\sigma(T)$ the spectrum of an operator $T\in \Bbb B(H)$.

\vspace{.2in}

\section{Results}
\setcounter{equation}{0}

%We begin with the following elementary lemma and give the proof for the reader's convenience. In the rest of the paper we will use the term oblique projection to denote a projection which need not be orthogonal.

The following elementary lemma is probably well-known, but we nevertheless give the proof for the reader's convenience. Here and in the rest of the paper we will use the term oblique projection to denote a projection which need not be orthogonal.
%
%This is probably a known result, but we are not aware of its existence. We give the proof for the reader's convenience.
%
%To our knowledge,

\begin{lemma}\label{elementarna}
Let $H$ and $K$ be Hilbert spaces. Suppose that $T\in \Bbb B(H,K)$ and $S\in \Bbb B(K,H)$ satisfy $ST=I$. Then
\begin{enumerate}
\item[(i)]
$
\mbox{Ker}\,S=(I-TS)(\mbox{Ker}\,T^*),
$
\item[(ii)]
$K=\mbox{Im}\,T \stackrel{.}{+}\mbox{Ker}\,S$,
\item[(iii)] $TS$ is the oblique projection onto $\mbox{Im}\,T$ parallel to $\mbox{Ker}\,S$.
\end{enumerate}
\end{lemma}
\proof
Let us first prove
\begin{equation}\label{prvaizleme}
\mbox{Ker}\,S=\mbox{Im}\,(I-TS).
\end{equation}
Indeed, from $ST=I$ we get $STS=S$ and $S(I-TS)=0$. This immediately implies $\mbox{Im}\,(I-TS) \subseteq \mbox{Ker}\,S$. Conversely, for $y \in \mbox{Ker}\,S$ we have $(I-TS)y=y$ which gives us $y \in \mbox{Im}\,(I-TS)$.

Next we claim
\begin{equation}\label{trecaizleme}
\mbox{Im}\,(I-TS)=(I-TS)(\mbox{Ker}\,T^*).
\end{equation}

Observe that (i) follows directly from (\ref{prvaizleme}) and (\ref{trecaizleme}).

To prove (\ref{trecaizleme}), we first note that the assumed equality $ST=I$ implies that $T$ has closed range. This may be argued by noting that $T$ is bounded from bellow. Alternatively, a more direct argument is as follows: if $(Tx_n)$ is a convergent sequence, it suffices to show that $(x_n)$ converges and this is immediate from $S(Tx_n)=x_n$.

Now take any $(I-TS)y\in \mbox{Im}\,(I-TS)$. Since the range of $T$ is closed, we may write $y=Tx+z$ for some $Tx\in \mbox{Im}\, T$ and $z\in
\mbox{Ker}\,T^*$. Next we observe that $(I-TS)(Tx)=Tx-T(ST)x=Tx-Tx=0$. From this we conclude  $(I-TS)y=(I-TS)(Tx+z)=(I-TS)z$. Thus,
$\mbox{Im}\,(I-TS)\subseteq (I-TS)(\mbox{Ker}\,T^*)$. Since the reverse inclusion is obvious, this completes the proof of (\ref{trecaizleme}) and hence of (i).

Let us prove (ii). Take any $y\in \mbox{Im}\,T\cap \mbox{Ker}\,S$. This means that $y=Tx$ for some $x$ and $Sy=0$. Thus, $x=STx=Sy=0$ which implies
$y=Tx=0$. Next, take arbitrary $y\in K$. As in the first part of the proof we have $y=Tx+z$ with some $Tx\in \mbox{Im}\,T$, $z\in \mbox{Ker}\,T^*$ and \begin{equation}\label{petanova}
(I-TS)y=(I-TS)z.
\end{equation}
Put $u=TSy\in \mbox{Im}\,T$ and $v=(I-TS)z$. Using (i), we have $v=(I-TS)z\in (I-TS)(\mbox{Ker}\,T^*)=\mbox{Ker}\,S$. Hence, we may rewrite (\ref{petanova}) in the form
$$
y=TSy+(I-TS)z=u+v\in \mbox{Im}\,T \stackrel{.}{+}\mbox{Ker}\,S
$$
which completes the proof of (ii).

To prove (iii), first note that $(TS)^2=TSTS=T(ST)S=TS$ which shows that $TS$ is an oblique projection. Obviously, $TS$ acts as the identity on $\mbox{Im}\,T$ and trivially on $\mbox{Ker}\,S$.
\qed

\vspace{.1in}

We are now ready for our main results.

\vspace{.1in}

\begin{theorem}\label{isti visak}
Let $(f_n)$ and $(g_n)$ be frames in $H$ dual to each other. Then $e(f_n)=e(g_n)$.
\end{theorem}
\proof
Let us denote the analysis operators of  $(f_n)$ and $(g_n)$ by $U$ and $V$, respectively. We must prove that
$
\mbox{dim}(\mbox{Ker}\,U^*)=\mbox{dim}(\mbox{Ker}\,V^*)$.
By the assumed duality we have $V^*U=I$.
Lemma \ref{elementarna} (i) (with $S=V^*$ and $T=U$) gives us
\begin{equation}\label{jezgre}
\mbox{Ker}\,V^* =(I-UV^*)(\mbox{Ker}\,U^*).
\end{equation}
This implies
$$\mbox{dim}(\mbox{Ker}\,V^*)=\mbox{dim}((I-UV^*)(\mbox{Ker}\,U^*))\leq \mbox{dim}(\mbox{Ker}\,U^*).$$

The opposite inequality follows by symmetry, since $V^*U=I$ is equivalent to $U^*V=I$.
\qed

\vspace{.2in}

\begin{cor}\label{dualsynthesis}
Let $(f_n)$ and $(g_n)$ be frames in $H$ dual to each other. Let $U$ and $V$ denote the analysis operators of $(f_n)$ and $(g_n)$, respectively. Then
\begin{itemize}
\item[(a)] $l^2=\mbox{Im}\,U \stackrel{.}{+}\mbox{Ker}\,V^*$,
\item[(b)] $UV^*$ is the oblique projection onto $\mbox{Im}\,U$ parallel to $\mbox{Ker}\,V^*$,
\item[(c)] $l^2=\mbox{Im}\,V \stackrel{.}{+}\mbox{Ker}\,U^*$,
\item[(d)] $VU^*$ is the oblique projection onto $\mbox{Im}\,V$ parallel to $\mbox{Ker}\,U^*$.
\end{itemize}
\end{cor}
\proof
(a) and (b) are immediate from Lemma \ref{elementarna} (ii) and (iii) with $S=V^*$ and $T=U$.
Since $V^*U=I$ is equivalent to $U^*V=I$, (c) and (d) follow from (a) and (b) by symmetry.
\qed

\vspace{.1in}

One should observe that, if we take the canonical dual $((U^*U)^{-1}f_n)$ of $(f_n)$ and its analysis operator $V=U(U^*U)^{-1}$, then the projection $UV^*$ from the above statement (b) becomes $UV^*=U(U^*U)^{-1}U^*$ which is precisely the orthogonal projection onto $\mbox{Im}\,U$.

As a consequence of the preceding discussion,
we obtain for an arbitrary frame $(f_n)$ in $H$ a general form of the synthesis operator of any frame that is dual to $(f_n)$. Denote by $U$ the analysis operator of $(f_n)$. Then, by Proposition 5.6.4 from \cite{Chr},  $V\in \Bbb B(H,l^2)$ is the analysis operator of a frame $(g_n)$ dual to $(f_n)$ if and only if its adjoint $V^*$ ({\it i.e.}, the corresponding synthesis operator) is given by $V^*=(U^*U)^{-1}U^*+W^*Q$, where
$Q\in \Bbb B(l^2)$ is the orthogonal projection to $(\mbox{Im}\,U)^{\perp}$ and $W\in \Bbb B(H,l^2)$ is arbitrary. (It is useful to note that $U(U^*U)^{-1}$ is in fact the pseudo-inverse $(U^*)^{\dag}$ of $U^*$.)

In our next corollary we provide another general form of $V^*$.

\begin{cor}\label{dualform}
Let $(f_n)$ be a frame in $H$ with the analysis operator $U$.
Then $V\in \Bbb B(H,l^2)$ is the analysis operator of a frame $(g_n)$ in $H$ dual to $(f_n)$ if and only if its adjoint $V^*$ is given by
$V^*=(U^*U)^{-1}U^*F$, where $F\in \Bbb B(l^2)$ is an oblique projection onto  $\mbox{Im}\,U$.
\end{cor}
\proof
Consider $V\in \Bbb B(H,l^2)$ whose adjoint operator $V^*$ is given by $V^*=(U^*U)^{-1}U^*F$, where $F\in \Bbb B(l^2)$ is an oblique projection onto  $\mbox{Im}\,U$.
Then $FU=U$ and hence $V^*U=(U^*U)^{-1}U^*FU=(U^*U)^{-1}U^*U=I$.

Conversely, suppose that $V\in \Bbb B(H,l^2)$ has the property $V^*U=I$. Then we know from the preceding corollary that
$l^2=\mbox{Im}\,U \stackrel{.}{+}\mbox{Ker}\,V^*$ and that $UV^*$ is the oblique projection onto $\mbox{Im}\,U$ parallel to $\mbox{Ker}\,V^*$.
Put $F=UV^*$. Then $(U^*U)^{-1}U^*F=(U^*U)^{-1}U^*UV^*=V^*$.

Observe that in the statement of the corollary and in the first part of the proof we have tacitly assumed that $l^2=\mbox{Im}\,U \stackrel{.}{+}X$ and that $F$ projects onto $\mbox{Im}\,U$ parallel to $X$, for some direct complement $X$ of $\mbox{Im}\,U$. The subspace $X$ played no other role in the proof. However, $X$ has to be closed in $l^2$ in order to ensure boundedness of $F$.

\qed

\vspace{.1in}

\begin{remark}
{\em Let $(f_n)$ be an arbitrary frame in $H$ with the analysis operator $U$. The preceding two corrolaries establish a $1-1$ correspondence between its dual frames and bounded oblique projections to $\mbox{Im}\,U$ ({\em i.e.} closed direct complements of $\mbox{Im}\,U$ in $l^2$).}
\end{remark}

\vspace{.1in}

\begin{remark}
{\em Theorem \ref{isti visak} is also true for finite frames in finite-dimensional Hilbert spaces. The same proof applies after observing that the analysis operators $U$ and $V$ take values in $\Bbb C^m$ (or $\Bbb R^m$), provided that frames under consideration are sequences with $m$ elements.}
\end{remark}

\vspace{.1in}

If we have a frame $(f_n)$  in $H$ and a surjection $T\in \Bbb B(H,K)$, then it is well known that
$(Tf_n)$ is a frame in $K$. Observe that for each $y\in K$ and $n\in \Bbb N$ we have $\langle y,Tf_n\rangle=\langle T^*y,f_n\rangle$. From this we conclude: if $U$ denotes the analysis operator of the original frame $(f_n)$, then the analysis operator $V$ of
$(Tf_n)$ is given by $V=UT^*$. Thus, $V^*=TU^*$ which obviously implies $\mbox{Ker}\,U^*\subseteq  \mbox{Ker}\,V^*$
and hence $e(f_n)\leq e(Tf_n)$.

We say that frames $(f_n)$ and $(g_n)$ for $H$ and $K$, respectively, are {\it equivalent} if there exists an invertible operator $T\in \Bbb B(H,K)$ such that $g_n=Tf_n,\,\forall n\in \Bbb N$. Since we then also have $f_n=T^{-1}g_n,\,\forall n\in \Bbb N$, the preceding discussion shows that equivalent frames have the same excess.

\vspace{.1in}

Recall now the concept of pseudo-duality for frames. Suppose that $(f_n)$ and $(g_n)$ are frames in $H$ with the analysis operators $U$ and $V$, respectively. We say that $(f_n)$ and $(g_n)$ are {\it pseudo-dual} to each other if $V^*U$ is an invertible operator. Note that $V^*Ux=\sum_{n=1}^{\infty}\langle x,f_n\rangle g_n,\,\forall x \in H$. This gives us, for each $x\in H$,
$$
x=(V^*U)((V^*U)^{-1}x)=\sum_{n=1}^{\infty}\langle (V^*U)^{-1}x,f_n\rangle g_n
=\sum_{n=1}^{\infty}\langle x,(U^*V)^{-1}f_n\rangle g_n.
$$
This shows that frames $((U^*V)^{-1}f_n)$ and $(g_n)$ are dual to each other. Since $(f_n)$ and
$((U^*V)^{-1}f_n)$ are equivalent, the preceding discussion and Theorem \ref{isti visak} give us $e(g_n)=e((U^*V)^{-1}f_n)=e(f_n)$.

Finally, recall from \cite{CrLa} that frames $(f_n)$ and $(g_n)$ in $H$ are called {\it approximately dual} to each other if the corresponding analysis operators $U$ and $V$ satisfy
$\|V^*U-I\|<1$. Clearly, approximate duality of frames implies pseudo-duality.

After all, we have proved the following corollary.

\vspace{.1in}

\begin{cor}\label{pseudo}
Pseudo-dual frames have the same excess. In particular, approximately dual frames have the same excess.
\end{cor}

\vspace{.1in}

Let us now make a comment on frame perturbations.

\begin{remark}
{\em
There is a number of results on frame perturbations. Typically, such results tell us that a sequence which is in a certain sense close to a frame must be a frame itself, in most cases sharing the same properties with the given one. We refer the reader to Chapter 15 in \cite{Chr} for a nice exposition of that part of the theory of frames.

Here we mention Theorem 1 in \cite{Chr1} and Theorem 2 in \cite{CasChr2} as important examples of perturbation results. Basically, the proofs of these results consist of two steps.
First, one shows that the sequence $(g_n)$ obtained by perturbing a frame $(f_n)$ is Bessel, so that its analysis operator $V$ is well defined and bounded and, secondly, that $V^*U$ is an invertible operator (where $U$ is the analysis operator of $(f_n)$).
From that, one easily concludes that  $(g_n)$ is indeed a frame.

Since $V^*U$ is invertible, this frame is pseudo-dual to $(f_n)$.
By the preceding corollary, this implies that $(f_n)$ and $(g_n)$ have the same excess. This last conclusion is already known in both cases; see Theorem 3.6 in \cite{CasChr} and Theorem 2 in \cite{CasChr2}, but the preceding discussion provides another (common) viewpoint for understanding the reason for the equality of excesses in these and similar situations.
}
\end{remark}

\vspace{.1in}

We now turn to another question in frame theory where the excess of a frame under consideration plays an important role. Suppose we are given a frame
$(f_n)$ in $H$. One may ask if there exists a Parseval dual for $(f_n)$; {\it i.e.}, a Parseval frame $(g_n)$ in $H$ such that $x=\sum_{n=1}^{\infty}\la x,f_n\ra g_n,\,\,\forall x \in H$. A necessary and sufficient condition is obtained in \cite{H}: $(f_n)$ does possess a Parseval dual if and only if $(f_n)$ can be obtained by applying an oblique projection to an orthonormal basis of a larger Hilbert space $K$ that contains $H$ as a subspace. Conveniently enough, this property of frames is cha\-rac\-te\-rized in \cite{ACRS}. So, by combining the results from these two papers, one obtains as a corollary the desired description of frames possessing Parseval duals. This is already noted in \cite{H}. It turns out that the key property is "sufficiently large excess". For reader's convenience here we state the characterization theorem and provide another, more direct proof.

\begin{theorem}\label{argentinski}
Let $(f_n)$ be a frame in $H$ with frame bounds $A$ and $B$ and the analysis operator $U$. Then $(f_n)$ possesses a Parseval dual if and only if the following two conditions are satisfied:
\begin{itemize}
\item[(a)] $A \geq 1$,
\item[(b)] $\mbox{dim}(\mbox{Im}\,(U^*U-I))\leq e(f_n)$.
\end{itemize}
\end{theorem}
\proof
Suppose first that $(g_n)$ is a Parseval dual for $(f_n)$; denote its analysis operator by $V$. As we have already mentioned in the discussion preceding Corollary \ref{dualform}, $V$ must be of the form $V=U(U^*U)^{-1}+QW$, where $Q\in \Bbb B(l^2)$ is the orthogonal projection to $(\mbox{Im}\,U)^{\perp}$ and $W\in \Bbb B(H,l^2)$ is arbitrary.

Since $(g_n)$ is a Parseval frame, we have $V^*V=I$ {\it i.e.},
$$((U^*U)^{-1}U^*+W^*Q)(U(U^*U)^{-1}+QW)=I.$$
Since $QU=0$ and $U^*Q=0$, this gives us
\begin{equation}\label{duplow}
(U^*U)^{-1} + W^*QW=I.
\end{equation}
In particular, this implies $(U^*U)^{-1}\leq I$. Therefore, the upper bound of the spectrum of $(U^*U)^{-1}$ is less than or equal to $1$; in other words, $\frac{1}{A}\leq 1$. Thus, we have proved $A \geq 1$.

Consider the orthogonal decomposition $H=\mbox{Ker}\,(U^*U-I) \oplus \overline{\mbox{Im}\,(U^*U-I)}$ and write, accordingly, $U^*U=I\oplus T$ (where $T\in \Bbb B(\overline{\mbox{Im}\,(U^*U-I)})$ denotes the operator induced by $U^*U$ on the invariant subspace $\overline{\mbox{Im}\,(U^*U-I)}$).
Observe that $T\geq 0$ and $\sigma(T) \subseteq [1,B]$.

We also claim that $I-T^{-1}$ is an injection. To see this, assume that $(I-T^{-1})x=0$ for some $x \in \overline{\mbox{Im}\,(U^*U-I)}$. This means that $T^{-1}x=x$ which implies $Tx=x$; in other words,
$U^*Ux=x$. Hence $x\in\mbox{Ker}\,(U^*U-I)$. Thus, $x=0$.

By (\ref{duplow}), we have $W^*QW=I-(U^*U)^{-1}=0\oplus (I-T^{-1})$. Observe that $\sigma(T^{-1})\subseteq [\frac{1}{B},1]$. We now claim that the range of $I-T^{-1}$ is a dense subspace of $\overline{\mbox{Im}\,(U^*U-I)}$. This is obvious if $I-T^{-1}$ is invertible. If $I-T^{-1}$ is not an invertible operator then, since it is self-adjoint and, by the preceding paragraph, injective, $-1$ has to be in the continuous part of its spectrum. Thus, we have proved that
the range of $W^*QW=I-(U^*U)^{-1}=0\oplus (I-T^{-1})$ is dense in $\overline{\mbox{Im}\,(U^*U-I)}$. This gives us
\begin{align*}
\mbox{dim}(\mbox{Im}\,(U^*U-I)) &= \overline{\mbox{dim}(\mbox{Im}\,(U^*U-I))} = \mbox{dim}(\mbox{Im}\,W^*QW)\leq \mbox{dim}(\mbox{Im}\,Q) \\
&= \mbox{dim}((\mbox{Im}\,U)^{\perp}) = e(f_n).
\end{align*}
To prove the converse, assume (a) and (b). Write again $U^*U=I\oplus T$ according to the decomposition $H=\mbox{Ker}\,(U^*U-I) \oplus \overline{\mbox{Im}\,(U^*U-I)}$. Again, $T\geq 0$ and, since $A \geq 1$, here we also have $\sigma(T) \subseteq [1,B]$.

Consider a continuous function $g:[1,\infty) \rightarrow [0,1)$ defined by $g(t)=\sqrt{1-\frac{1}{t}}$. Put $G=g(T)$. We now use the assumption (b) to find a partial isometry $R\in \Bbb B(H,l^2)$ whose initial space is $\overline{\mbox{Im}\,(U^*U-I)}$ with final space {\it contained} in $\mbox{Ker}\,U^*=
(\mbox{Im}\,U)^{\perp}$. Finally, denote by $P \in \Bbb B(H)$ the orthogonal projection onto $\overline{\mbox{Im}\,(U^*U-I)}$.

Let $V=U(U^*U)^{-1}+R(0\oplus G)P$. Then
\begin{align*}
V^*V &= \left((U^*U)^{-1}U^*+P(0\oplus G)R^*\right)\left(U(U^*U)^{-1}+R(0\oplus G)P\right) \\
&= (U^*U)^{-1}+(0\oplus G^2)=(I\oplus T^{-1})+(0\oplus (I-T^{-1})) = I.
\end{align*}
Let us now put $g_n=V^*e_n,\,n\in \Bbb N$, where $(e_n)$ is the canonical orthonormal basis in $l^2$.
Since $V^*V=I$, the sequence $(g_n)$ is a Parseval frame in $H$. Obviously, we also have $V^*U=I$ which means that $(g_n)$ is a dual frame to $(f_n)$.
\qed

\vspace{.1in}

The above condition (a) is not crucial, since it can be ensured by rescaling the original frame (although, the construction then yields only a tight dual frame with the frame bound different from $1$.) Condition (b) is essential; it tells us that the excess should be at least large as $d=\mbox{dim}(\mbox{Im}\,(U^*U-I))=\mbox{dim}((\mbox{Ker}\,(U^*U-I))^{\perp})$. Note that the number $d$ can be interpreted as a kind of a measure of deviation of the original frame from being Parseval. Namely, the characterizing Parseval property - $U^*U=I$ - is trivially fulfilled on the subspace $\mbox{Ker}\,(U^*U-I)$. So, any deviation from the Parseval property has its origin in the orthogonal complement $(\mbox{Ker}\,(U^*U-I))^{\perp}$.

We also refer the reader to \cite{BB} for a discussion of some related properties of frames with $d<\infty$.

\vspace{.1in}

Next we turn to the so called fundamental identity for Parseval frames. Let $(f_n)$ be an arbitrary Parseval frame in $H$. It is proved in \cite{BCEK} that for each $J\subset \Bbb N$ and for all $x\in H$ we have
\begin{equation}\label{fundamental}
\sum_{n\in J}|\langle x,f_n\rangle|^2+\left\|\sum_{n\in \Bbb N\setminus J}\langle x,f_n\rangle f_n\right\|^2 =
\sum_{n\in \Bbb N \setminus J}|\langle x,f_n\rangle|^2+\left\|\sum_{n\in J}\langle x,f_n\rangle f_n\right\|^2
\end{equation}
and
\begin{equation}\label{tricetvrtine}
\frac{3}{4}\|x\|^2\leq \sum_{n\in J}|\langle x,f_n\rangle|^2+\left\|\sum_{n\in \Bbb N\setminus J}\langle x,f_n\rangle f_n\right\|^2
\leq \|x\|^2.
\end{equation}
Given a Parseval frame $(f_n)$ and $J\subseteq \Bbb N$, one can define the quantities
$$
\nu_-(J)=\mbox{inf}\left\{\frac{1}{\|x\|^2}\left(\sum_{n\in J}|\langle x,f_n\rangle|^2+\left\|\sum_{n\in \Bbb N\setminus J}\langle x,f_n\rangle f_n\right\|^2\right): x \not =0\right\}
$$
and
$$
\nu_+(J)=\mbox{sup}\left\{\frac{1}{\|x\|^2}\left(\sum_{n\in J}|\langle x,f_n\rangle|^2+\left\|\sum_{n\in \Bbb N\setminus J}\langle x,f_n\rangle f_n\right\|^2\right): x \not =0\right\}.
$$
Clearly, (\ref{tricetvrtine}) implies
\begin{equation}\label{ni}
\frac{3}{4}\leq \nu_-(J)\leq\nu_+(J)\leq 1,\,\forall J\subset \Bbb N.
\end{equation}
These inequalities are discussed in \cite{G}. Here we provide a related property of frames with finite excess. Loosely speaking, in Theorem \ref{ocjena} and Corollary \ref{posljedica} below we prove: if $(f_n)$ is a Parseval frame with finite excess then for each $\epsilon >0$ we have $\nu_-(J)\geq(1-\epsilon)$ for the majority of subsets $J$ of $\Bbb N$.

\vspace{.1in}

\begin{theorem}\label{ocjena}
Let $(f_n)$ be a Parseval frame in  $H$ such that $e(f_n)<\infty$. Then for each $\epsilon >0$ there exists $n_0 \in \Bbb N$ such that
$$
\sum_{n=1}^{n_0}|\langle x,f_n\rangle|^2+\left\|\sum_{n=n_0+1}^{\infty}\langle x,f_n\rangle f_n\right\|^2 > (1-\epsilon)\|x\|^2,\,\,\forall x \in H,
$$
and
$$
\sum_{n=n_0+1}^{\infty}|\langle x,f_n\rangle|^2+\left\|\sum_{n=1}^{n_0}\langle x,f_n\rangle f_n\right\|^2> (1-\epsilon)\|x\|^2,\,\,\forall x \in H.
$$
\end{theorem}
\proof

Let $U$ be the analysis operator of the frame $(f_n)$ and let $(e_n)$ be the canonical orthonormal basis for $l^2$. Put $\mbox{Im}\,U=M$ and denote by $P\in \Bbb B(l^2)$ the orthogonal projection onto $M$. Clearly, $(Pe_n)$ is a Parseval frame for $M$, unitarily equivalent to $(f_n)$. Hence, without loss of generality, we may assume that $H=M$ and $f_n=Pe_n,\,\forall n \in \Bbb N$.

Fix $\epsilon >0$. Recall from Proposition 5.5 in \cite{CHet_al} that $\sum_{n=1}^{\infty}(1-\|f_n\|^2)=e(x_n)$.
Since by hypothesis we have $e(x_n)=k<\infty$, there exists $n_0 \in \Bbb N$
such that
\begin{equation}\label{jedan}
\sum_{n=n_0+1}^{\infty}(1-\|f_n\|^2)<\epsilon.
\end{equation}
Let $\{w_1,w_2,\ldots, w_k\}$ be an orthonormal basis for $M^{\perp}$. Then, for each $x \in l^2$, we have
$$
\left\|(I-P)\left(\sum_{n=n_0+1}^{\infty}\langle x,e_n\rangle e_n\right)\right\|^2=\sum_{j=1}^k\left|\left\langle \sum_{n=n_0+1}^{\infty}\langle x,e_n\rangle e_n, w_j\right\rangle\right|^2=
$$
$$
\sum_{j=1}^k\left|\sum_{n=n_0+1}^{\infty}\langle x,e_n\rangle\langle e_n,w_j\rangle\right|^2\leq
\sum_{j=1}^k\left(\sum_{n=n_0+1}^{\infty}|\langle x,e_n\rangle|^2\right)\left(\sum_{n=n_0+1}^{\infty}|\langle e_n, w_j\rangle|^2\right)\leq
$$
$$
\|x\|^2\sum_{j=1}^k\left(\sum_{n=n_0+1}^{\infty}|\langle e_n, w_j\rangle|^2\right)=\|x\|^2\sum_{n=n_0+1}^{\infty}\left(\sum_{j=1}^k
\langle e_n, w_j\rangle|^2
\right)=
$$
$$
\|x\|^2\sum_{n=n_0+1}^{\infty}\|(I-P)e_n\|^2=\|x\|^2\sum_{n=n_0+1}^{\infty}(1-\|Pe_n\|^2)=
$$
$$
\|x\|^2\sum_{n=n_0+1}^{\infty}(1-\|f_n\|^2)
\stackrel{(\ref{jedan})}{<}\epsilon \|x\|^2.
$$
Thus, we have proved
\begin{equation}\label{dva}
\left\|(I-P)\left(\sum_{n=n_0+1}^{\infty}\langle x,e_n\rangle e_n\right)\right\|^2<\epsilon \|x\|^2,\,\,\forall x \in \l^2.
\end{equation}
Let us now take any $x\in M$. Observe that $Px=x$. Then we have
$$
\sum_{n=1}^{n_0}|\langle x,f_n\rangle|^2+\left\|\sum_{n=n_0+1}^{\infty}\langle x,f_n\rangle f_n\right\|^2=
$$
$$
\sum_{n=1}^{n_0}|\langle x,Pe_n\rangle|^2+\left\|\sum_{n=n_0+1}^{\infty}\langle x,Pe_n\rangle Pe_n\right\|^2=
$$
$$
\sum_{n=1}^{n_0}|\langle Px,e_n\rangle|^2+\left\|P\left(\sum_{n=n_0+1}^{\infty}\langle Px,e_n\rangle e_n\right)\right\|^2=
$$
$$
\sum_{n=1}^{n_0}|\langle x,e_n\rangle|^2+\left\|P\left(\sum_{n=n_0+1}^{\infty}\langle x,e_n\rangle e_n\right)\right\|^2=
$$
$$
\sum_{n=1}^{n_0}|\langle x,e_n\rangle|^2+\left\|\sum_{n=n_0+1}^{\infty}\langle x,e_n\rangle e_n\right\|^2-
\left\|(I-P)\left(\sum_{n=n_0+1}^{\infty}\langle x,e_n\rangle e_n\right)\right\|^2=
$$
$$
\|x\|^2-\left\|(I-P)\left(\sum_{n=n_0+1}^{\infty}\langle x,e_n\rangle e_n\right)\right\|^2\stackrel{(\ref{dva})}{>}\|x\|^2-\epsilon \|x\|^2
=(1-\epsilon)\|x\|^2.
$$
This proves the first inequality. The second inequality now follows from (\ref{fundamental}).

\qed

\vspace{.1in}

Notice that inequality (\ref{dva}) in the above proof remains true if we replace the summation over the set $\{n_0+1,n_0+2,\ldots \}$ by the sumation over any subset of $\Bbb N$ disjoint from $\{1,2,\ldots ,n_0\}$ - the same proof applies. Thus, we have the following corollary:

\begin{cor}\label{posljedica}
Let $(f_n)$ be a Parseval frame for $H$ such that $e(f_n)<\infty$. Then for each $\epsilon >0$ there exists $n_0 \in \Bbb N$ such that for each subset $J$ of $\Bbb N$ with the property $\{1,2,\ldots ,n_0\}\subseteq J$ we have
$$
\sum_{n\in J}|\langle x,f_n\rangle|^2+\left\|\sum_{n\in \Bbb N\setminus J}\langle x,f_n\rangle f_n\right\|^2 > (1-\epsilon)\|x\|^2,\,\,\forall x \in H,
$$
and
$$
\sum_{n\in \Bbb N \setminus J}|\langle x,f_n\rangle|^2+\left\|\sum_{n\in J}\langle x,f_n\rangle f_n\right\|^2> (1-\epsilon)\|x\|^2,\,\,\forall x \in H.
$$
\end{cor}

\vspace{.1in}

\begin{ex}
{\em
Let $a=(\alpha_n)\in l^2$ such that $\|a\|^2=\sum_{n=1}^{\infty}|\alpha_n|^2=1$, assume additionally that $|\alpha_1|^2> \frac{7}{8}$ and $\alpha_n\not =0,\,\forall n \in \Bbb N$.

Let $M=\{x\in l^2:\langle x,a\rangle=0\}=\{a\}^{\perp}$. Denote by $P\in \Bbb B(l^2)$ the orthogonal projection onto $M$. Let $(e_n)$ be the canonical orthonormal basis for $l^2$. Put $x_n=Pe_n,\,n\in \Bbb N$. Then $(x_n)$ is a Parseval frame for $M$. Its analysis operator is the inclusion map $Ux=x,\,x\in M$; thus, the synthesis operator $U^*$ coincides with $P$ regarded as a map from $l^2$ onto $M$. From $\mbox{Ker}\,P=\mbox{span}\{a\}$ we conclude $e(x_n)=1$. (Note in passing that here one obtains a Riesz basis by removing just one element from $(x_n)$. In fact, since we assumed $\alpha_n\not =0,\,\forall n \in \Bbb N$, any $x_n$ can be removed. We omit the details.)

Observe that $$\|x_n\|^2=\|Pe_n\|^2=1-\|(I-P)e_n\|^2=1-|\langle e_n,a\rangle|^2=1-|\alpha_n|^2,\,\forall n \in \Bbb N.$$
Thus, by Proposition 5.5 from \cite{CHet_al}, we have
$$
1=e(x_n)=\sum_{n=1}^{\infty}(1-\|x_n\|^2)=\sum_{n=1}^{\infty}|\alpha_n|^2.
$$
This gives us
$$
\sum_{n=2}^{\infty}(1-\|x_n\|^2)=\sum_{n=2}^{\infty}|\alpha_n|^2=1-|\alpha_1|^2<\frac{1}{8}.
$$
Thus, keeping notations from the proof of Theorem \ref{ocjena}, here we have, for $\epsilon =\frac{1}{8}$, $n_0=1$. Since for each subset $J$ of $\Bbb N$ we have either $1\in J$ or $1\in \Bbb N \setminus J$, Corollary \ref{posljedica} implies
$$
\sum_{n\in J}|\langle x,x_n\rangle|^2+\left\|\sum_{n\in \Bbb N\setminus J}\langle x,x_n\rangle x_n\right\|^2 > (1-\epsilon)\|x\|^2
=\frac{7}{8}\|x\|^2,\,\,\forall x \in H,\,\forall J\subseteq \Bbb N.
$$
In particular, this implies
$$
\nu_-:=\mbox{inf}\,\left\{\nu_-(J):J\subseteq \Bbb N\right\}\geq \frac{7}{8}.
$$
Notice that, by (\ref{ni}), we have for each Parseval frame $\nu_-\geq \frac{3}{4}$. Our frame $(x_n)$ constructed here serves as an example of a Parseval frame for which $\nu_-$ is strictly greater than $\frac{3}{4}$. A characterization of Parseval frames with the property $\nu_->\frac{3}{4}$ is not known.
}
\end{ex}

\end{document}